\documentclass[11pt]{amsart}
\usepackage{amsmath, amsfonts, amsbsy, amssymb, upref, chatterjee1}

\newcommand{\mr}{\mathcal{R}}

\begin{document}
\title{The Ghirlanda-Guerra identities without averaging}
\author{Sourav Chatterjee}
\date{November 18, 2009}
\thanks{The author's research was supported in part by NSF grant DMS 0707054 and a Sloan Research Fellowship}
\address{Courant Institute of Mathematical Sciences, New York University, 251 Mercer Street, New York, NY 10012}
\begin{abstract}
The Ghirlanda-Guerra identities are one of the most mysterious features of spin glasses. We prove the GG identities in a large class of models that includes the Edwards-Anderson model, the random field Ising model, and the Sherrington-Kirkpatrick model in the presence of a random external field.  Previously, the GG identities were rigorously proved only  `on average' over a range of temperatures or under small perturbations. 
\end{abstract}
\maketitle

\section{Introduction}
Consider the Sherrington-Kirpatrick model of spin glasses \cite{sk75} in the presence of an external field, defined as follows. Let $N$ be a positive integer. Let $(g_{ij})_{1\le i<j\le N}$ and $(g_i)_{1\le i\le N}$ be collections of i.i.d.\ standard Gaussian random variables. Fix three real numbers $\beta, \gamma, h$. Define a Hamiltonian $H_N$ on $\{-1,1\}^N$ as 
\[
-H_N(\bos) := \frac{\beta}{\sqrt{N}} \sum_{1\le i<j\le N} g_{ij} \sigma_i \sigma_j + h \sum_{1\le i\le N} \sigma_i + \gamma \sum_{1\le i\le N} g_i \sigma_i. 
\]
Let $G_N$ be the (random) probability measure on $\{-1,1\}^N$ that puts mass 
\[
G_N(\{\bos\}) \propto e^{-H_N(\bos)}
\]
at each configuration $\bos$. This is the Gibbs measure of the SK model. As usual, for a function $f:(\{-1,1\}^N)^n \ra \rr$ we write
\[
\smallavg{f} := \sum_{\bos^1,\ldots,\bos^n} f(\bos^1,\ldots,\bos^n) G_N(\{\bos^1\})\cdots G_N(\{\bos^n\}). 
\]
Following the notation in Talagrand~\cite{talagrand03}, we write
\[
\nu(f) := \ee\smallavg{f}. 
\]
Let $\bos^1,\bos^2,\ldots$ be i.i.d.\ configurations from the Gibbs measure. In the spin glass parlance, these are known as replicas. The overlap between two replicas $\bos^\ell$ and $\bos^{\ell'}$ is defined as
\[
R_{\ell,\ell'} := \frac{1}{N}\sum_{i=1}^N \sigma_i^\ell \sigma_i^{\ell'}. 
\]
Now fix a number $n$, and a bounded measurable function $f:\rr^{n(n-1)/2}\ra \rr$. Let $\mr_n$ denote the collection $(R_{\ell,\ell'})_{1\le \ell<\ell'\le n}$. Writing $f = f(\mr_n)$ for simplicity, let
\begin{equation}\label{deltadef}
\delta_N(\beta,\gamma, h) := \nu(R_{1,n+1} f) - \frac{1}{n}\nu(R_{1,2}) \nu(f) - \frac{1}{n}\sum_{2\le \ell\le n} \nu(R_{1,\ell} f).
\end{equation}
The Ghirlanda-Guerra identities \cite{ghirlandaguerra98} claim that for any such $n$ and $f$, for almost all $\beta, \gamma, h$, 
\begin{equation}\label{ggid}
\lim_{N\ra \infty} \delta_N(\beta, \gamma, h) = 0. 
\end{equation}
There is an extended form of the GG identities that claims more: in the infinite volume limit, the conditional distribution of $R_{1,n+1}$ given $\mr_n$ is 
\[
\frac{1}{n}\ml_{R_{1,2}} + \frac{1}{n}\sum_{2\le \ell\le n} \delta_{R_{1,\ell}},
\]
where $\ml_{R_{1,2}}$ is the (unconditional) law of $R_{1,2}$, and $\delta_{R_{1,\ell}}$ is the point mass at $R_{1,\ell}$. Clearly, \eqref{ggid} is a special case of this extended version. 

There is an early version of the GG identities that are sometimes called the Aizenman-Contucci identities; these appeared in \cite{aizenmancontucci98}. The paper \cite{aizenmancontucci98} also introduced the important notion of stochastic stability in spin glass models. An interesting connection between stochastic stability and the GG identities was discovered by Arguin \cite{arguin08}. 

The importance of the GG identities stems from their universal nature. The identities are supposed to hold in a wide array of spin glass models, both mean field and short range. A striking recent discovery of Panchenko~\cite{panchenko09}, following on the insightful work of Arguin and Aizenman \cite{arguinaizenman09}, indicates that the GG identities may actually solve the long-standing mystery of ultrametricity in spin glasses. 

It is not difficult to show (see e.g.\ Talagrand \cite{talagrand03}, Section 2.12) that for any $a > 0$, and any $\beta, h$, 
\begin{equation}\label{av}
\lim_{N\ra\infty} \int_{-a}^a |\delta_N(\beta, x, h)| dx = 0.
\end{equation}
Similar results were proved for a general class of spin glass models, including the short range Edwards-Anderson model, by Contucci and Giardin\`a~\cite{contuccigiardina07}. However, this does not prove \eqref{ggid}.  The best advance till date is due to Talagrand~\cite{talagrand09} who showed that for any $\beta, \gamma, h$, there exists a sequence $\gamma_N \ra \gamma$ such that
\[
\lim_{N\ra \infty} \delta_N(\beta, \gamma_N, h) = 0. 
\]
Let $\psi_N$ be the free energy per particle, defined as
\begin{equation*}\label{free}
\psi_N := \frac{1}{N}\log \sum_{\bos} e^{-H_N(\bos)}. 
\end{equation*}
Let $p_N = p_N(\beta, \gamma, h) := \ee(\psi_N)$. As easy extension of the Guerra-Toninelli argument \cite{guerratoninelli02} shows that $\lim_{N\ra \infty} p_N(\beta, \gamma, h)$ exists for all $\beta, \gamma, h$. Let us call this limit $p(\beta, \gamma, h)$. Then  $p$ is a convex function of each argument. 
\begin{thm}\label{thm1}
Suppose $p$ is differentiable in $\gamma$ at a certain value of $(\beta, \gamma, h)$, where $\gamma \ne 0$. Then at this value the Ghirlanda-Guerra identities hold, that is, $\lim_{N\ra\infty} \delta_N(\beta, \gamma, h) = 0$. 
\end{thm}
Clearly, the major drawback of the above result is that it only works for $\gamma \ne 0$. This is problematic since $\gamma = 0$ is the most interesting case. A minor drawback is that the theorem holds for almost every value of $\gamma$ instead of all values. But this is not a serious problem since it is not clear whether the GG identities hold at points of phase transition, and the original paper of Ghirlanda and Guerra \cite {ghirlandaguerra98} does not make such a claim either. And lastly, it will be an important breakthrough if one can prove a similar statement for the extended GG identities, but that seems to be out of reach at present.

The proof of Theorem~\ref{thm1} admits a vast generalization, which we state below. 

For each $N$, suppose $\mu_N$ is a probability measure on $\{-1,1\}^N$. Suppose $\ma_N$ is a sequence of finite sets such that $|\ma_N| = O(N)$. For each $N$, suppose $(f_\alpha)_{\alpha \in \ma_N}$ is a collection of functions from $\{-1,1\}^N$ into $[-1,1]$ and $(g_{\alpha})_{\alpha\in \ma_N}$ is a collection of i.i.d.\ standard Gaussian random variables. Fix a parameter~$\gamma$. Let $G_N$ be the probability measure on $\{-1,1\}^N$ satisfying
\[
G_N(\{\bos\}) \propto \mu_N(\{\bos\}) \exp\biggl(\gamma\sum_{\alpha\in \ma_N} g_\alpha f_\alpha(\bos)\biggr). 
\]
Let $\psi_N$ be the free energy per particle in this model, that is,
\begin{equation}\label{psidef}
\psi_N = \psi_N(\gamma) := \frac{1}{N}\log \sum_{\bos}\mu_N(\bos)\exp\biggl(\gamma\sum_{\alpha\in \ma_N} g_\alpha f_\alpha(\bos)\biggr). 
\end{equation}
Let $p_N = \ee(\psi_N)$. Assume that $p(\gamma) := \lim_{N\ra \infty} p_N(\gamma)$ exists. 

Next, for two replicas $\bos^1$ and $\bos^2$, define the generalized overlap as 
\[
R_{1,2} := \frac{1}{N}\sum_{\alpha\in \ma_N} f_\alpha(\bos^1) f_\alpha(\bos^2).
\]
Assume that $R_{1,1}$ is a deterministic constant for all $\bos^1$. 
Given $n$ and a bounded measurable function $f:\rr^{n(n-1)/2}\ra \rr$, let $\delta_N(\gamma)$ be defined as in~\eqref{deltadef}, in terms of the generalized overlaps defined above. 
\begin{thm}\label{thm2}
Suppose $p$ is differentiable at a point $\gamma$. Then the Ghirlanda-Guerra identities hold at this $\gamma$, that is,  $\lim_{N\ra\infty} \delta_N(\gamma)=0$. 
\end{thm}
Note that this result includes the Sherrington-Kirkpatrick model and Derrida's $p$-spin models under a random external field. It also covers lattice models like the random field Ising model and the Edwards-Anderson model~\cite{edwardsanderson75}. For example, in the Edwards-Anderson model, the spins are located on a subset of the $d$-dimensional integer lattice; $\ma_N$ is the set of edges in this subset, and for an edge $\alpha = (i,j)\in \ma_N$, $f_\alpha(\bos)= \sigma_i \sigma_j$. The measure $\mu_N$ is just the uniform distribution on $\{-1,1\}^N$ in this case. In this model, our generalized overlap is simply what's known as the bond overlap, and $\gamma$ is the inverse temperature parameter. Here we do not need to apply a random external field. The existence of the limit of the free energy per particle in the lattice models follow from simple arguments, e.g.\ those in~\cite{ruelle99}, Chapter~2. 

\section{Proof}
We will only prove Theorem \ref{thm2}, since Theorem \ref{thm1} is a special case of Theorem \ref{thm2}. Let
\[
H = H(N) := \frac{\sum g_\alpha f_\alpha(\bos)}{N}.
\]
It is easy to prove via integration by parts (see e.g.\ \cite{talagrand03}, Section~2.12) that the Ghirlanda-Guerra identities hold if 
\begin{equation}\label{s0}
\lim_{N\ra\infty} \nu(|H-\nu(H)|) = 0. 
\end{equation}
We prove \eqref{s0} in two steps. First, we  show that 
\begin{equation}\label{s1}
\lim_{N\ra\infty} \ee|\smallavg{H}- \nu(H)| = 0. 
\end{equation}
This part requires no new ideas. The proof goes as follows. Note that
\begin{equation}\label{diff}
\smallavg{H} = \psi_N'(\gamma) \ \text{ and } \  \nu(H) = p_N'(\gamma).  
\end{equation}
Fix some $\gamma' > \gamma$. By the convexity of $\psi_N$ and \eqref{diff}, we have
\begin{equation}\label{upbd}
\smallavg{H} \le \frac{\psi_N(\gamma')-\psi_N(\gamma)}{\gamma'-\gamma}. 
\end{equation}
Again, by standard concentration of measure (see e.g.\ \cite{talagrand03}, Theorem 2.2.4) we know that for any $\gamma$, 
\begin{equation}\label{conc}
\ee|\psi_N(\gamma) - p_N(\gamma)| \le \sqrt{\var(\psi_N(\gamma))}\le  C|\gamma| \sqrt{\frac{|\ma_N|}{N^2}},
\end{equation}
where $C$ is a universal constant. Thus, 
\begin{align*}
&\ee\biggl| \frac{\psi_N(\gamma')-\psi_N(\gamma)}{\gamma'-\gamma} - p'(\gamma)\biggr| \\
&\le \frac{C(|\gamma| + |\gamma'|)}{\gamma'-\gamma}\sqrt{\frac{|\ma_N|}{N^2}}  + \frac{|p_N(\gamma') - p(\gamma')| + |p_N(\gamma) - p(\gamma)|}{\gamma'-\gamma} \\
&\qquad + \biggl|\frac{p(\gamma') - p(\gamma)}{\gamma'-\gamma} - p'(\gamma)\biggr|. 
\end{align*}
Since $|\ma_N|=O(N)$ and $p_N\ra p$ pointwise, we get
\begin{align*}
\limsup_{N\ra \infty} \ee\biggl| \frac{\psi_N(\gamma')-\psi_N(\gamma)}{\gamma'-\gamma} - p'(\gamma)\biggr|&\le \biggl|\frac{p(\gamma') - p(\gamma)}{\gamma'-\gamma} - p'(\gamma)\biggr|.
\end{align*}
Combining with \eqref{upbd}, we see that 
\begin{align*}
\limsup_{N\ra \infty} \ee(\smallavg{H}- p'(\gamma))_+ &\le \biggl|\frac{p(\gamma') - p(\gamma)}{\gamma'-\gamma} - p'(\gamma)\biggr|.
\end{align*}
Here $x_+$ denotes the positive part of a real number $x$. Since this bound holds for any $\gamma' > \gamma$ and $p$ is differentiable at $\gamma$, we get
\[
\lim_{N\ra \infty} \ee(\smallavg{H}  - p'(\gamma))_+ = 0. 
\]
Similarly, considering $\gamma' < \gamma$ and repeating the steps, we can show that the limit of the negative part is zero as well. Thus,
\[
\lim_{N\ra \infty} \ee|\smallavg{H}  - p'(\gamma)| = 0. 
\]
By Jensen's inequality, this gives
\[
\lim_{N\ra \infty} |\nu(H)- p'(\gamma)| = 0. 
\]
Combining the last two identities, we get \eqref{s1}. This completes the first step. In the next part of the proof, we show that
\begin{equation}\label{s2}
\lim_{N\ra\infty} \nu(|H- \smallavg{H}|) = 0. 
\end{equation}
Combined with \eqref{s1}, this will complete the proof of \eqref{s0} and hence the theorem. To show \eqref{s2}, it suffice to prove that 
\begin{equation}\label{toprove}
\lim_{N\ra \infty} \ee\smallavg{(H-\smallavg{H})^2} = 0. 
\end{equation}
The key new idea at this stage is the use of the following result. Essentially, this is the Parseval identity for the $L^2$ norm of a Gaussian functional expressed as a sum of squares of its Fourier coefficients in the orthogonal basis of multidimensional Hermite polynomials.
\begin{thm}[\cite{chatterjee09}]\label{varformula}
Let $\bg = (g_1,\ldots, g_n)$ be a vector of i.i.d.\ standard Gaussian random variables, and let $f$ be a $C^\infty$ function of $\bg$ with bounded derivatives of all orders. Then 
\[
\var(f) = \sum_{k=1}^\infty \frac{1}{k!}\sum_{1\le i_1,\ldots,i_k\le n} \biggl(\ee\biggl(\frac{\partial^k f}{\partial g_{i_1}\cdots \partial g_{i_k}}\biggr)\biggr)^2. 
\]
The convergence of the infinite series is part of the conclusion. 
\end{thm}
The formula in its above form is stated and proved in \cite{chatterjee09}. Different versions of the identity were observed by a number of authors; let us refer to Section 3.8 of \cite{chatterjee09} for a discussion of these observations and alternate proofs.

We will now prove \eqref{toprove} using the above theorem. In the following, we will simply write $f_\alpha$ instead of $f_\alpha(\bos)$. Integration by parts gives
\begin{align*}
&\ee\smallavg{(H-\smallavg{H})^2} = \ee\smallavg{H^2} - \ee(\smallavg{H}^2)\\
&= \frac{1}{N^2}\sum_{\alpha, \alpha'} \ee\bigl(g_\alpha g_{\alpha'} \bigl(\smallavg{f_\alpha f_{\alpha'}} - \smallavg{f_\alpha}\smallavg{f_{\alpha'}}\bigr)\bigr)\\
&= \frac{1}{N^2}\sum_{\alpha,\alpha'} \ee\biggl(\mpar{}{g_\alpha}{g_{\alpha'}} \bigl(\smallavg{f_\alpha f_{\alpha'}} - \smallavg{f_\alpha}\smallavg{f_{\alpha'}}\bigr)\biggr) \\
&\qquad + \frac{1}{N^2}\sum_{\alpha} \ee(\smallavg{f_\alpha^2}-\smallavg{f_\alpha}^2). 
\end{align*}
Now note that
\[
\frac{1}{N}\bigl(\smallavg{f_\alpha f_{\alpha'}} - \smallavg{f_\alpha}\smallavg{f_{\alpha'}}\bigr) = \frac{1}{\gamma^2}\mpar{\psi_N}{g_\alpha}{g_{\alpha'}},
\]
where $\psi_N$ is the free energy per particle defined in \eqref{psidef}. Thus, by the Cauchy-Schwarz inequality and Theorem~\ref{varformula}, we get
\begin{align*}
\ee\smallavg{(H-\smallavg{H})^2} &\le \frac{1}{N\gamma^2}\sum_{\alpha, \alpha'} \ee\biggl(\frac{\partial^4 \psi_N}{\partial g_\alpha^2 \partial g_{\alpha'}^2}\biggr) + \frac{2|\ma_N|}{N^2}\\
&\le \frac{|\ma_N|}{N\gamma^2 } \biggl(\sum_{\alpha, \alpha'} \biggl(\ee\biggl(\frac{\partial^4 \psi_N}{\partial g_\alpha^2 \partial g_{\alpha'}^2}\biggr)\biggr)^2\biggr)^{1/2} + O(1/N)\\
&\le \frac{|\ma_N|\sqrt{4! \var(\psi_N)}}{N\gamma^2} + O(1/N). 
\end{align*}
The proof is now completed by applying the variance  inequality \eqref{conc}.
\vskip.2in
\noindent{\bf Acknowledgments.} The authors thanks Michel Talagrand and Pierluigi Contucci for helpful comments.

\end{document}